\documentclass{article}
\usepackage[T1]{fontenc}
\usepackage{geometry}
\usepackage{amssymb}
\usepackage{amsmath,cases}
\usepackage{amsthm}
\usepackage{mathrsfs}
\usepackage{enumitem}
\usepackage{xcolor}
\usepackage{graphicx}
\usepackage{psfrag}
\usepackage{textcomp}
\usepackage{caption}
\usepackage{url}
\usepackage{calc}
\usepackage{bbm}
\usepackage[active]{srcltx}   
\usepackage{pdfsync}
\usepackage{todonotes}
\usepackage[colorlinks=true,bookmarksnumbered=true, pdfauthor={J\'er\'emie Bettinelli, Eric Fusy, C\'ecile Mailler, Lucas Randazzo}]{hyperref}									

\newcommand{\eps}{\varepsilon}
\newcommand{\de}{\mathrel{\mathop:}\hspace*{-.6pt}=}

\newcommand{\bC}{\mathbb{C}}
\newcommand{\bD}{\mathbb{D}}
\newcommand{\dsk}{\overline{\mathbb{D}}}
\newcommand{\be}{\mathbbm{e}}
\newcommand{\bR}{\mathbb{R}}

\newcommand{\cB}{\mathcal{B}}
\newcommand{\cK}{\mathcal{K}}
\newcommand{\cP}{\mathcal{P}}

\theoremstyle{plain}
\newtheorem{thm}{Theorem}
\newtheorem{alg}{Algorithm}

\newtheorem{prop}[thm]{Proposition}

\theoremstyle{definition}
\newtheorem*{rem}{Remark}

\newenvironment{pre}[1][\proofname]{%
  \proof[#1]%
}{\endproof}

\makeatletter
\renewcommand*{\@fnsymbol}[1]{\ensuremath{\ifcase#1\or \dagger\or \ddagger\or
   \mathsection\or \mathparagraph\or \|\or **\or \dagger\dagger
   \or \ddagger\ddagger \else\@ctrerr\fi}}
\makeatother
\title{Convergence of uniform noncrossing partitions toward the Brownian triangulation}
\author{J\'er\'emie Bettinelli\thanks{cnrs \& Laboratoire d'Informatique de l'\'Ecole polytechnique; \href{mailto:jeremie.bettinelli@normalesup.org}{\nolinkurl{jeremie.bettinelli@normalesup.org}};\newline \nolinkurl{www.normalesup.org/}\texttildelow\nolinkurl{bettinel}. This work is partially supported by Grant ANR-14-CE25-0014 (GRAAL).}}

\begin{document}
\maketitle

\begin{abstract}
We give a short proof that a uniform noncrossing partition of the regular $n$-gon weakly converges toward Aldous's Brownian triangulation of the disk, in the sense of the Hausdorff topology. This result was first obtained by Curien \& Kortchemski, using a more complicated encoding. Thanks to a result of Marchal on strong convergence of Dyck paths toward the Brownian excursion, we furthermore give an algorithm that allows to recursively construct a sequence of uniform noncrossing partitions for which the previous convergence holds almost surely.

In addition, we also treat the case of uniform noncrossing pair partitions of even-sided polygons.

\end{abstract}


\section{Introduction}\label{secintro}

Configurations of noncrossing diagonals of a regular polygon have been the focus of many studies, from a geometrical, from an enumerative and from a probabilistic point of view. Various natural models have been studied (see for instance~\cite{CuKo14} and the references therein). Among these models, 
noncrossing partitions are of particular interest as they bear many applications in a wide range of areas; see for instance~\cite{McC06} for a survey of the topic.

\begin{figure}[ht]
		\psfrag{0}[B][B]{$0$}
		\psfrag{1}[B][B]{$1$}
		\psfrag{2}[B][B]{$2$}
		\psfrag{3}[B][B]{$3$}
		\psfrag{4}[B][B]{$4$}
		\psfrag{5}[B][B]{$5$}
		\psfrag{6}[B][B]{$6$}
		\psfrag{7}[B][B]{$7$}
		\psfrag{8}[B][B]{$8$}
		\psfrag{9}[B][B]{$9$}
	\centering\includegraphics[width=.95\linewidth]{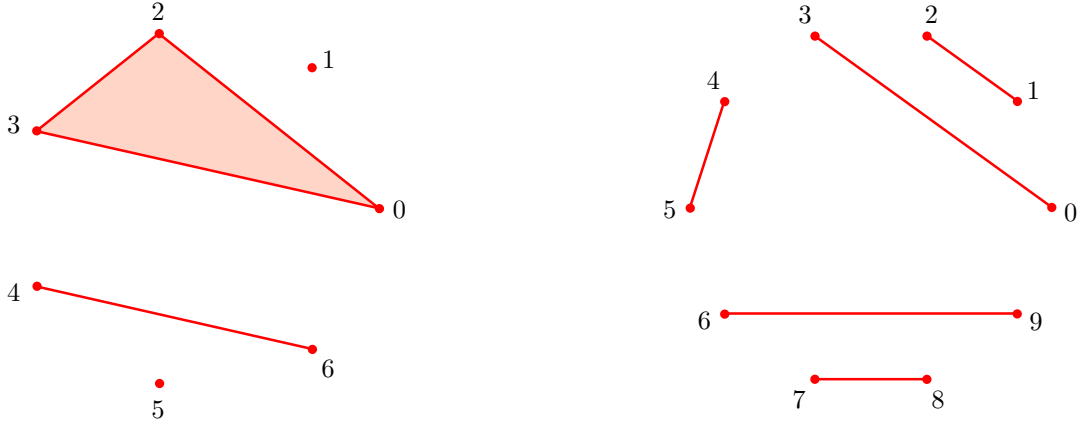}
	\caption{\textbf{Left.} The noncrossing partition $\big\{\{\omega_7^0,\omega_7^2,\omega_7^3\},\{\omega_7^1\},\{\omega_7^4,\omega_7^6\},\{\omega_7^5\}\big\}$ of~$P_7$. \textbf{Right.} The noncrossing pair partition $\big\{\{\omega_{10}^0,\omega_{10}^3\},\{\omega_{10}^1,\omega_{10}^2\},\{\omega_{10}^4,\omega_{10}^5\},\{\omega_{10}^6,\omega_{10}^9\},\allowbreak\{\omega_{10}^7,\omega_{10}^8\}\big\}$ of~$P_{10}$.}
	\label{ncp}
\end{figure}

We denote by~$P_n$ the regular $n$-gon of the complex plane with vertex coordinates $\omega_n^k\de e^{2\mathrm{i}\pi k/n}$, $k\in \{0,1,\ldots, n-1\}$. A \emph{noncrossing partition} of~$P_n$ is a partition of the set $\{\omega_n^0,\ldots,\omega_n^{n-1}\}$ such that the convex hulls of its blocks are pairwise disjoint; a \emph{noncrossing pair partition} of~$P_n$ is a noncrossing partition of~$P_n$ whose blocks are all of size exactly~$2$ (see Figure~\ref{ncp}). Note that the latter only exists for even values of~$n$.

A (geodesic) \emph{lamination} is a closed subset of the unit disk $\dsk \de \big\{ z \in \bC : |z| \leq 1\big\}$ that can be expressed as a union of chords whose intersections with the open unit disk~$\bD$ are pairwise disjoint. For two complex numbers~$z$, $z'$, we denote by $[z,z']$ the Euclidean line segment of the complex plane joining~$z$ with~$z'$. A noncrossing partition can thus be seen as a lamination as follows. With the block $\big\{\omega_n^{i_1},\omega_n^{i_2},\ldots,\omega_n^{i_k}\big\}$, where $i_1<i_2<\ldots<i_k$, we associate the polygon $[\omega_n^{i_1},\omega_n^{i_{2}}]\cup[\omega_n^{i_2},\omega_n^{i_{3}}]\cup\ldots\cup[\omega_n^{i_k},\omega_n^{i_{1}}]$. The lamination is then defined as the union over the partition blocks of the associated polygons.

Recall that the \emph{Hausdorff distance} between two closed subsets~$A$, $B \subseteq \dsk$ is defined as
\[
\inf\Big\{ \eps >0\,:\, A \subseteq B^{(\eps)} \text{ and } B \subseteq A^{(\eps)}\Big\},
\]
where, for any $X \subseteq \dsk$, we denoted by $X^{(\eps)}\de\{z\in\dsk\,:\, d(z,X)\le \eps\}$ the $\eps$-enlargement of~$X$. Endowed with the Hausdorff metric, the set of all closed subsets of~$\dsk$ is a compact metric space. Moreover, it is not hard to check that the set of all laminations is a closed, thus compact, subset of this metric space.

We are interested in the limit of sequences of larger and larger noncrossing partitions, seen as laminations, for the Hausdorff topology. The above setting was proposed by Aldous \cite{aldous94,aldous94b} for the study of random triangulations; it was later used by Kortchemski~\cite{Kor14stablelam}, Curien \& Kortchemski \cite{CuKo14}, Curien \& Le~Gall \cite{CuLG11tri}, Kortchemski \& Marzouk \cite{KoMa16,KoMa17} for the study of many models, including uniform noncrossing partitions and uniform noncrossing pair partitions. 

In many cases, the limiting object is a random compact set called the \emph{Brownian triangulation} $\cB$ defined as follows. Let $(\be_t)_{0\le t\le 1}$ be a normalized Brownian excursion\footnote{A \emph{normalized Brownian excursion} is a standard Brownian motion on $[0,1]$ starting from~$0$ and conditioned on being at~$0$ at time~$1$ and staying positive on $(0,1)$. As this is a zero-probability event, some care is needed for a proper definition: see e.g.\ \cite[Chapter~XII]{revuz99cma}.} and, for $s$, $t\in[0,1]$, declare $s \stackrel\be\sim  t$ whenever $\be_s=\be_t = \min_{\min(s,t)\le r\le  \max(s,t)} \be_r$. The
Brownian triangulation is the set
\begin{equation}\label{defB}
\cB \de \bigcup_{s \stackrel\be\sim  t} \big[e^{2\mathrm{i} \pi s},e^{2\mathrm{i} \pi t}\big]\,.
\end{equation}

The set~$\cB$ is almost surely a closed subset of~$\dsk$ and furthermore a continuous triangulation of~$\dsk$, in the sense that each connected component of $\dsk\setminus\cB$ is an open Euclidean triangle whose vertices belong to the unit circle~\cite{legall08slb}.

\begin{thm}[{\cite[Theorem~3.8]{CuKo14}}]\label{thmpn}
Let $\cP_n$ (resp.\ $\tilde\cP_n$) be a random variable uniformly distributed over the set of all noncrossing partitions of~$P_{n}$ (resp.\ noncrossing pair partitions of~$P_{2n}$), seen as a lamination. Then~$\cP_{n}$ and~$\tilde\cP_n$ both weakly converge toward the Brownian triangulation, for the Hausdorff topology.
\end{thm}

The point of the present work is to provide a more straightforward proof of the previous theorem. For this reason, we chose not to include too many historical references on the subject; we refer the reader to the above references and references therein for more details. In~\cite{CuKo14}, the authors first notice that a noncrossing partition of~$P_n$ is close to a noncrossing pair partition of~$P_{2n}$; as a consequence, the result for noncrossing pair partitions implies the result for noncrossing partitions. They then encode a noncrossing pair partition by a tree, which they further encode by a Dyck path; after proper scaling, this path converges to the normalized Brownian excursion. This approach is quite robust but needs in particular a nontrivial result stating that the leaves of a conditioned Galton--Watson tree are asymptotically uniformly spread on the tree.

Instead, we will give a direct encoding of a noncrossing partition of~$P_n$ by a Dyck path and conclude in a more straightforward manner. In addition, we will give below a recursive construction of noncrossing partitions and noncrossing pair partitions that converge almost surely. The algorithms we propose are the transcription in terms of partitions of an algorithm on Dyck paths due to Marchal~\cite{Marchal03}, which itself is the transcription in terms of Dyck paths of R{\'e}my's famous algorithm~\cite{remy85} on trees. Let us also mention at this point that Curien \& Le~Gall \cite{CuLG11tri} also study some sequences of laminations obtained by a recursive construction but, in their case, the limiting object is a continuous triangulation that differs from the Brownian triangulation.

\bigskip

\begin{figure}[ht]
		\psfrag{0}[B][B]{$0$}
		\psfrag{1}[B][B]{$1$}
		\psfrag{2}[B][B]{$2$}
		\psfrag{3}[B][B]{$3$}
		\psfrag{4}[B][B]{$4$}
		\psfrag{5}[B][B]{$5$}
		\psfrag{6}[B][B]{$6$}
		\psfrag{7}[B][B]{$7$}
		\psfrag{8}[B][B]{$8$}
		\psfrag{9}[B][B]{$9$}		
	\centering\includegraphics[width=.4\linewidth]{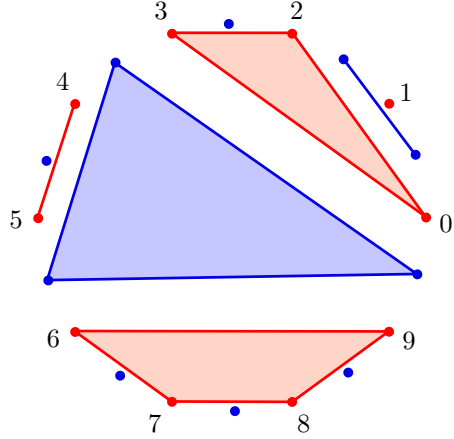}
	\caption{The Kreweras complement of the partition \textcolor{red}{$\big\{\{\omega_{10}^0,\omega_{10}^2,\omega_{10}^3\},\{\omega_{10}^1\},\allowbreak\{\omega_{10}^4,\omega_{10}^5\},\{\omega_{10}^6,\omega_{10}^7,\omega_{10}^8,\omega_{10}^9\}\big\}$} is \textcolor{blue}{$\big\{\{\omega_{20}^1,\omega_{20}^3\},\{\omega_{20}^5\},\{\omega_{20}^7,\omega_{20}^{11},\omega_{20}^{19}\},\{\omega_{20}^9\},\{\omega_{20}^{13}\},\allowbreak\{\omega_{20}^{15}\},\{\omega_{20}^{17}\}\big\}$}. For instance, the blocks \textcolor{red}{$\{\omega_{10}^0,\omega_{10}^2,\omega_{10}^3\}$} and \textcolor{blue}{$\{\omega_{20}^5\}$} are neighbors.}
	\label{compl}
\end{figure}

Before presenting our growing algorithms, let us define the \emph{Kreweras complement} of a noncrossing partition~$\cP$ of~$P_n$ as the partition~$\cK$ of the set $\{\omega_{2n}^1,\omega_{2n}^3,\ldots,\omega_{2n}^{2n-1}\}$ whose blocks are given by the connected components of the complement in~$\dsk$ of the lamination corresponding to~$\cP$ (see Figure~\ref{compl}). Up to rotation of $-\pi/n$, the Kreweras complement of a noncrossing partition of~$P_n$ is also a noncrossing partition of~$P_n$. Taking again the Kreweras complement of the latter noncrossing partition of~$P_n$ and rotating it by an angle of $\pi/n$ yields back the original noncrossing partition.

\begin{figure}[ht]
		\psfrag{1}[B][B][.65]{$1$}
		\psfrag{2}[B][B][.65]{$2$}
		\psfrag{3}[B][B][.65]{$3$}
		\psfrag{a}[B][B][.65]{$\ge 2$}
		\psfrag{b}[B][B][.65]{+$1$}
		\psfrag{k}[B][B]{$k$}
		\psfrag{l}[B][B]{$l$}
		\psfrag{i}[B][B][.8]{\itshape inserting a vertex}
		\psfrag{s}[B][B][.8]{\itshape slicing}
	\centering\includegraphics[width=.95\linewidth]{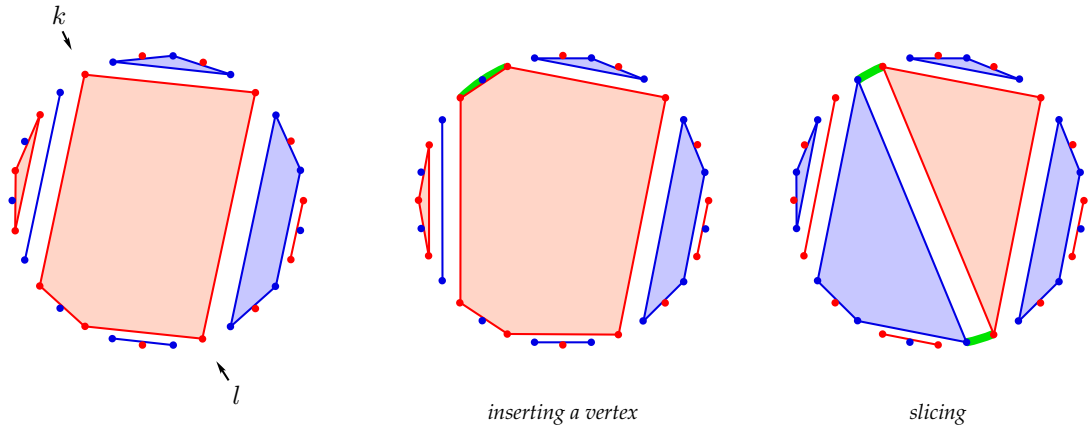}
	\caption{The noncrossing partitions obtained by inserting a vertex and by slicing at position~$k$.}
	\label{grow}
\end{figure}

We will see two possible ways of defining a noncrossing partition of~$P_{n+1}$ from a noncrossing partition~$\cP$ of~$P_n$ and an index $k\in\{0,1,\ldots,2n\}$ (see Figure~\ref{grow}). The first operation consists in adding two vertices between~$\omega_{2n}^k$ and~$\omega_{2n}^{k+1}$, and declaring the second vertex as belonging to the block of~$\omega_{2n}^k$. We then remap the $2n+2$ vertices onto the $2n+2$-th roots of unity in such a way that the cyclic order is preserved and~$\omega_{2n}^k$ is mapped to~$\omega_{2n+2}^k$. The resulting noncrossing partition of~$P_{n+1}$ is said to be obtained from~$\cP$ by \emph{inserting a vertex at position~$k$}.

For the second operation, we need to consider the last element in counterclockwise order before~$\omega_{2n}^{2n}$ that belongs to the same block of $\cP\cup\cK$ as~$\omega_{2n}^k$: let
\[
l\de \max\{j\le 2n\,:\,\omega_{2n}^{k} \text{ and } \omega_{2n}^{j}\text{ are in the same block of } \cP\cup\cK\}\,.
\]
We split each of the vertices~$\omega_{2n}^k$ and~$\omega_{2n}^l$ into two new vertices and remap the resulting $2n+2$ vertices onto the $2n+2$-th roots of unity in such a way that the cyclic order is preserved and~$\omega_{2n}^0$ is mapped to~$\omega_{2n+2}^0$. We define a noncrossing partition of~$P_{n+1}$ by declaring any two $n+1$-th roots of unity to be in the same block whenever their preimages were in the same block of $\cP\cup\cK$. Note that the Kreweras complement of this noncrossing partition is obtained in a similar manner by considering the other $2n+2$-th roots of unity. Our operation has the effect of slicing the block of~$\omega_{2n}^k$ along the chord $[\omega_{2n}^k,\omega_{2n}^l]$ into two blocks, one lying in the noncrossing partition and the other one lying in its Kreweras complement. We say that the resulting noncrossing partition of~$P_{n+1}$ is obtained from~$\cP$ by \emph{slicing at position~$k$}.

\begin{rem}
Observe that, whenever $l=k$, the noncrossing partitions obtained from~$\cP$ by slicing and by inserting a vertex at position~$k$ are the same.
Note also that the result is not the same when $k=0$ and when $k=2n$. In fact, the partitions obtained by insertion only differ by a rotation of $2\pi/n$ and the ones obtained by slicing differ by complement and rotation of~$\pi/n$.
\end{rem}

We iteratively construct a sequence of random noncrossing partitions~$(\cP_n)_{n\ge 1}$ using the following algorithm.
\begin{alg}[Constructing a sequence of noncrossing partitions]\label{algoncp}
\hspace{0pt}

\begin{enumerate}
	\item Let $\cP_1=\big\{\{\omega_1^0\}\big\}$ be the only partition of~$P_1$.
	\item Generate~$\cP_{n+1}$ from~$\cP_n$ as follows:
	\begin{enumerate}
		\item choose an integer~$k$ uniformly at random in $\{0,1,\ldots,2n\}$;
		\item with probabilities $1/2$ - $1/2$, set~$\cP_{n+1}$ to be obtained from~$\cP_n$ 
		\begin{itemize}
			\item either by inserting a vertex at position~$k$,
			\item or by slicing at position~$k$.
		\end{itemize}
	\end{enumerate}
\end{enumerate}
\end{alg}

\begin{prop}\label{propncp}
Let $(\cP_n)_{n\ge 1}$ be constructed by Algorithm~\ref{algoncp}. Then, for each~$n$, the partition~$\cP_n$ is uniformly distributed over the set of noncrossing partitions of~$P_n$. Moreover, seen as a lamination, $\cP_{n}$ almost surely converges toward the Brownian triangulation, for the Hausdorff topology.
\end{prop}

\begin{figure}[ht]
		\psfrag{1}[B][B][.65]{$1$}
		\psfrag{2}[B][B][.65]{$2$}
		\psfrag{3}[B][B][.65]{$3$}
		\psfrag{a}[B][B][.65]{$\ge 2$}
		\psfrag{b}[B][B][.65]{+$1$}
		\psfrag{k}[B][B]{$k$}
		\psfrag{l}[B][B]{$l$}
		\psfrag{i}[B][B][.8]{\itshape inserting a short chord}
		\psfrag{s}[B][B][.8]{\itshape inserting a long chord}
	\centering\includegraphics[width=.95\linewidth]{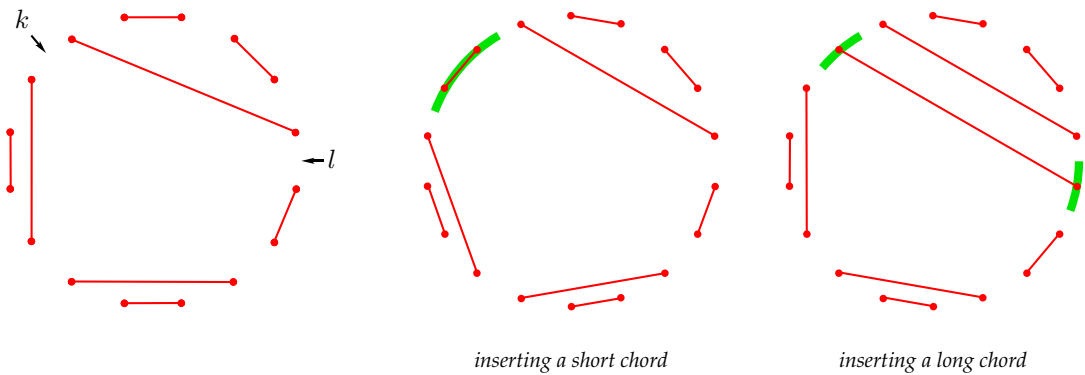}
	\caption{Growing a pair partition.}
	\label{growpair}
\end{figure}

We can play a similar game for noncrossing pair partitions (see Figure~\ref{growpair}). Let~$\cP$ be a noncrossing pair partition of~$P_{2n}$ and let $k\in\{0,1,\ldots,2n\}$. It will be more convenient to rotate the picture by an angle of~$-\pi/2n$, so that~$\cP$ is now a partition of $\big\{\omega_{4n}^j, j \text{ odd}\big\}$. We consider the Kreweras complement~$\cK$ of~$\cP$ and we set
\[
l\de \max\{j\le 2n\,:\,\omega_{2n}^{k} \text{ and } \omega_{2n}^{j}\text{ are in the same block of } \cK\}\,.
\]
We either add two vertices at the location of~$\omega_{2n}^{k}$ or one at the location of~$\omega_{2n}^{k}$ and one at the location of~$\omega_{2n}^{l}$. Then, in both cases, we declare the added vertices to form one new block and remap the $2n+2$ vertices onto the odd $4n+4$-th roots of unity in such a way that the cyclic order is preserved and~$\omega_{4n}^{2k-1}$ is mapped to~$\omega_{4n+4}^{2k-1}$. We say that the resulting noncrossing pair partitions of~$P_{2n+2}$ are obtained from~$\cP$ respectively by \emph{inserting a short chord} and by \emph{inserting a long chord} at position~$k$.

Similarly as above, we iteratively construct a sequence of random noncrossing pair partitions~$(\tilde\cP_n)_{n\ge 1}$.
\begin{alg}[Constructing a sequence of noncrossing pair partitions]\label{algopair}
\hspace{0pt}

\begin{enumerate}
	\item Let $\tilde\cP_1=\big\{\{\omega_2^0,\omega_2^1\}\big\}$ be the only pair partition of~$P_2$.
	\item Generate~$\tilde\cP_{n+1}$ from~$\tilde\cP_n$ as follows:
	\begin{enumerate}
		\item choose an integer~$k$ uniformly at random in $\{0,1,\ldots,2n\}$;
		\item with probabilities $1/2$ - $1/2$, set~$\tilde\cP_{n+1}$ to be obtained from~$\tilde\cP_n$ by inserting at position~$k$
		\begin{itemize}
			\item either a short chord,
			\item or a long chord.
		\end{itemize}
	\end{enumerate}
\end{enumerate}
\end{alg}

\begin{prop}\label{proppair}
Let $(\tilde\cP_n)_{n\ge 1}$ be constructed by Algorithm~\ref{algopair}. Then, for each~$n$, the partition~$\tilde\cP_n$ is uniformly distributed over the set of noncrossing pair partitions of~$P_{2n}$. Moreover, seen as a lamination, $\tilde\cP_{n}$ almost surely converges toward the Brownian triangulation, for the Hausdorff topology.
\end{prop}

The remainder of the paper is organized as follows. In Section~\ref{secenc}, we show how to encode a noncrossing partition by a Dyck path. Section~\ref{seccvd} is devoted to the proof of Theorem~\ref{thmpn} and Section~\ref{secas} to the proofs of Propositions~\ref{propncp} and~\ref{proppair}.

\section{Encoding noncrossing partitions by Dyck paths}\label{secenc}

\begin{figure}[ht]
		\psfrag{0}[B][B]{$0$}
		\psfrag{1}[B][B]{$0$}
		\psfrag{2}[B][B]{$1$}
		\psfrag{3}[B][B]{$2$}
		\psfrag{4}[B][B]{$3$}
	\centering\includegraphics[width=.95\linewidth]{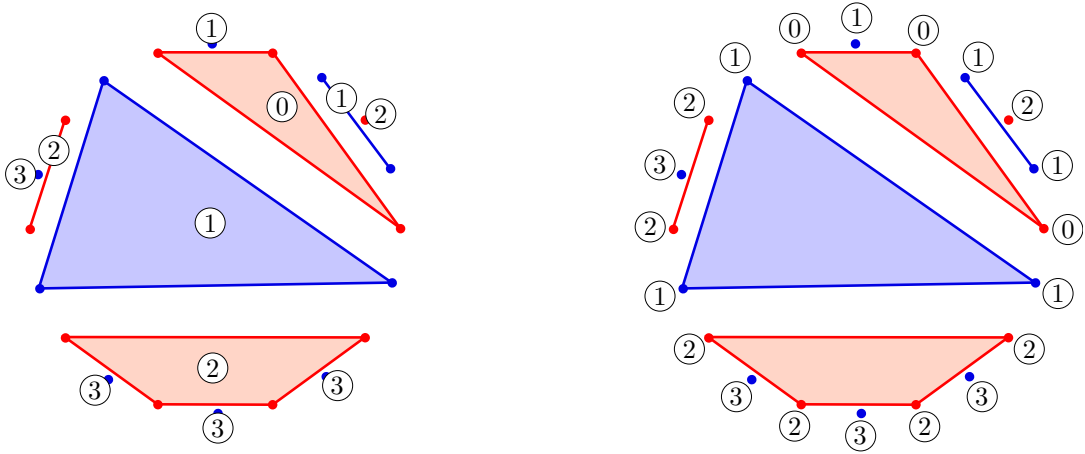}
	\caption{Encoding a noncrossing partition by a Dyck path.}
	\label{enc}
\end{figure}

We encode a noncrossing partition~$\cP$ of~$P_n$ by assigning integer labels to $\omega_{2n}^k$, $0\le k\le 2n-1$, as follows (see Figure~\ref{enc}). We let~$\cK$ be the Kreweras complement of~$\cP$ and we say that two blocks of $\cP\cup\cK$ are \emph{neighbors} if there exists an integer~$k$ such that~$\omega_{2n}^k$ belongs to one block and~$\omega_{2n}^{k+1}$ belongs to the other block. We first label the blocks of $\cP\cup\cK$ by assigning label~$0$ to the block that contains~$\omega_{2n}^0$ and, inductively, assigning label $\ell+1$ to each not yet labeled neighbor of a block labeled~$\ell$. We then assign to each~$\omega_{2n}^k$, $0\le k\le 2n-1$, the label of the block to which it belongs.

If we denote by~$\ell_k$ the label assigned to~$\omega_{2n}^k$, $0\le k \le 2n$, then the path $(\ell_0, \ell_1,\ldots,\ell_{2n})$ is a $2n$-step Dyck path\footnote{Recall that a Dyck path is a finite sequence $(a_0,a_1,\ldots,a_l)$ of nonnegative integers such that $a_0=a_l=0$ and $|a_{k+1}-a_k|=1$ for all $k\in\{0,1,\ldots,l-1\}$.}. Moreover, this operation yields a bijection between noncrossing partitions of~$P_n$ and $2n$-step Dyck paths; the inverse operation goes as follows. Let $(\ell_0, \ell_1,\ldots,\ell_{2n})$ be a $2n$-step Dyck path. Then the noncrossing partition is given by the equivalence classes of the relation
\begin{equation}\label{idn}
\omega_{2n}^i \sim \omega_{2n}^j \iff \ell_{i}=\ell_{j}=\min_{\min(i,j)\le k\le  \max(i,j)} \ell_{k}
\end{equation}
for even~$i$, $j$. Furthermore, the Kreweras complement of this noncrossing partition is given by the equivalence classes of~\eqref{idn} for odd~$i$, $j$. Note also that, if over the set $[\min(i,j), \max(i,j)]$, the function~$\ell$ only reaches its minimum at the extremities, then the chord $[\omega_{2n}^i,\omega_{2n}^j]$ belongs to the lamination.

\begin{rem}
Let us note that this encoding appeared in~\cite{Stump13}. It is also easy to check that this is in fact exactly the Dyck path that encodes the dual tree of the associated noncrossing pair partition (see Figure~\ref{bijpair} and~\cite{CuKo14}, in particular Figure~7) or, equivalently, the (properly rooted) dual two-type tree of~\cite{KoMa17} (see Figure~2 therein).
\end{rem}

\section{Convergence in distribution}\label{seccvd}

Let us start with uniform noncrossing partitions. The proof is very similar to that of~\cite{CuKo14} but circumvent some technicalities because of the encoding we use; we give it in full detail for the sake of self-containment.

\begin{pre}[Proof of Theorem~\ref{thmpn} for uniform noncrossing partitions]
For each $n\ge 1$, let~$\cP_n$ be a random variable uniformly distributed over the set of noncrossing partitions of~$P_n$ and let $L_n:[0,1]\to\bR_+$ be the function defined as follows. We consider the labeling of the $2n$-th roots of unity given by the encoding of Section~\ref{secenc}. For $k\in\{0,1,\ldots,2n\}$, we let $L_n(k/2n)$ be the label of~$\omega_{2n}^{k}$ and we extend~$L_n$ to $[0,1]$ by linear interpolation between these values.

A well-known conditioned version of Donsker's invariance principle due to Kaigh~\cite[Theorem~2.6]{kaigh76ipr} states that the following convergence holds in distribution for the uniform topology on the space of continuous real-valued functions on~$[0,1]$:
\begin{equation}\label{kaigh}
\left( \frac {L_n(s)}{\sqrt{2n}} \right)_{0 \le s \le 1} \to ( \be_s )_{0 \le s \le 1}\,.
\end{equation}

Using Skorokhod's representation theorem, we may and will assume that the previous convergence holds almost surely. As the set of all closed subsets of~$\dsk$ endowed with the Hausdorff metric is a compact metric space, it suffices to show that any accumulation point of $(\cP_n)_n$ is the Brownian triangulation~$\cB$, defined by~\eqref{defB}. Let~$\cP$ be such an accumulation point. 

We first claim that $\cB\subseteq \cP$ almost surely. It is a classical fact that the local minimums of~$\be$ on $(0,1)$ are almost surely distinct. On the set of full measure where this property holds, if $s\stackrel\be\sim t$ with $s<t$, we can always find even $s_n$, $t_n\in\{0,2,4,\ldots,2n\}$ such that $s_n<t_n$,
\[
\frac{s_n}{2n}\to s\,,\qquad \frac{t_n}{2n}\to t\qquad \text{and}\qquad L_n\Big(\frac{s_n}{2n}\Big)=L_n\Big(\frac{t_n}{2n}\Big)<\min_{[\frac{s_n+1}{2n},\frac{t_n-1}{2n}]} L_n\,.
\]
Indeed, this fact is plain if $\be_r > \be_s=\be_t$ for all $r \in (s,t)$; otherwise, there exists $r \in (s,t)$ such that $\be_r =
\be_s=\be_t$, so that~$\be_r$ is a local minimum and, as a result, neither~$s$ nor~$t$ are times of a local minimum and one can find $s'<s$ and $t'>t$ arbitrarily close to~$s$ and~$t$ such that $(s',t')$ falls into the previous case.
The latter property implies that the chord $\big[\omega_{2n}^{s_n},\omega_{2n}^{t_n}\big]$ belongs to~$\cP_n$ (recall~\eqref{idn} and the discussion thereafter). As a result, the segment $[e^{2\mathrm{i}\pi s},e^{2\mathrm{i}\pi t}]\subseteq\cP$ and the claim follows.  

Now, in order to see that $\cP\subseteq\cB$, observe that, as all the~$\cP_n$ are laminations, $\cP$ is also a lamination. Then, as~$\cB$ is almost surely maximal for the inclusion relation (\cite[Proposition~2.1]{legall08slb}), we necessarily have $\cB=\cP$. This completes the proof.
\end{pre}

\begin{figure}[ht]	
		\psfrag{0}[B][B]{$0$}
		\psfrag{1}[B][B]{$1$}
		\psfrag{2}[B][B]{$2$}
		\psfrag{3}[B][B]{$3$}
		\psfrag{4}[B][B]{$4$}
		\psfrag{5}[B][B]{$5$}
		\psfrag{6}[B][B]{$6$}
		\psfrag{7}[B][B]{$7$}
		\psfrag{8}[B][B]{$8$}
		\psfrag{9}[B][B]{$9$}		
	\centering\includegraphics[width=.95\linewidth]{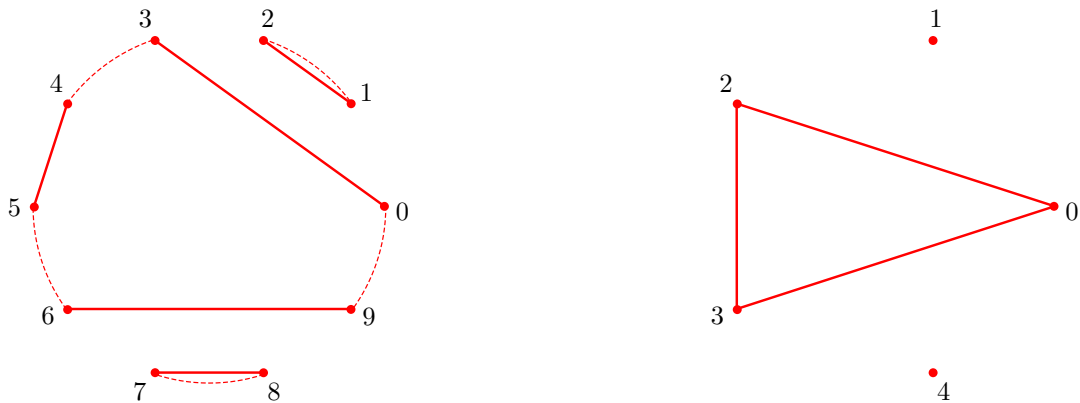}
	\caption{The noncrossing partition of~$P_{n}$ that corresponds to a noncrossing pair partition of~$P_{2n}$ is obtained by identifying $\omega_{2n}^{2k-1}$ with $\omega_{2n}^{2k}$ for each $1\le k\le n$.}
	\label{bijpair}
\end{figure}

We now turn to noncrossing pair partitions. We can use the observation of~\cite{CuKo14} that there exists a simple bijection between noncrossing partitions of~$P_n$ and noncrossing pair partitions of~$P_{2n}$, such that the corresponding laminations are at Hausdorff distance less than~$\pi/2n$ (see Figure~\ref{bijpair}).

Alternatively, we will see that the encoding of Section~\ref{secenc} is well behaved with respect to the property of being a noncrossing pair partition. 

\begin{pre}[Proof of Theorem~\ref{thmpn} for uniform noncrossing pair partitions]
Let~$\cP$ be a noncrossing partition of~$P_{2n}$. As in the introduction, we rotate the picture by an angle of~$-\pi/2n$ and see~$\cP$ as a partition of $\big\{\omega_{4n}^j, j \text{ odd}\big\}$ (plainly, this will bear no effects in the limit). We consider the Kreweras complement~$\cK$ of~$\cP$ and let $(\ell_0, \ell_1,\ldots,\ell_{4n})$ be its encoding Dyck path. By definition, $\cP$ is a noncrossing \emph{pair} partition if and only if the equivalence classes given by~\eqref{idn} for odd indices are all of size~$2$ (where~$n$ is replaced by~$2n$). Equivalently, for each $k\in\{0,1,\ldots,2n-1\}$, we have $|\ell_{2k+2}-\ell_{2k}|=2$, so that the path $(\ell_0/2,\ell_2/2,\ell_4/2,\ldots,\allowbreak\ell_{4n}/2)$ is a Dyck path that encodes all the information of the original Dyck path, and which does not satisfy any constraints (see Figure~\ref{encpair}).

\begin{figure}[ht]
		\psfrag{0}[B][B][.8]{$0$}
		\psfrag{1}[B][B][.8]{$1$}
		\psfrag{2}[B][B][.8]{$2$}
		\psfrag{3}[B][B][.8]{$3$}
	\centering\includegraphics[width=.4\linewidth]{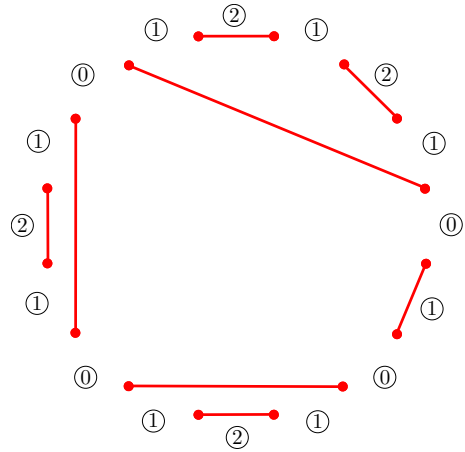}
	\caption{The unconstrained Dyck path that encodes a noncrossing pair partition.}
	\label{encpair}
\end{figure}

As a result, a uniform noncrossing pair partition of~$P_{2n}$ is encoded by a uniform $2n$-step Dyck path. By Kaigh's theorem, after proper rescaling, this path converges to the normalized Brownian excursion, and the original encoding $4n$-step Dyck path clearly converges toward the same limit, multiplied by~$\sqrt 2$. One then concludes exactly as in the case of uniform noncrossing partitions (the extra multiplicative factor does not alter the identifications).
\end{pre}

\section{Almost sure convergence}\label{secas}

In this section, we prove Propositions~\ref{propncp} and~\ref{proppair}. In fact, they are both straightforward consequences of~\cite{Marchal03} and the encoding we use. 

\begin{pre}[Proof of Proposition~\ref{propncp}]
In terms of encoding Dyck paths, inserting a vertex at position~$k$ amounts to inserting one up-step followed by one down-step right after time~$k$ and slicing at position~$k$ amounts to lifting up by one the part of the path between time~$k$ and the first subsequent time the path becomes strictly lower than its height at time~$k$ (see Figure~\ref{grow_l}).

\begin{figure}[ht]
		\psfrag{1}[B][B][.65]{$1$}
		\psfrag{2}[B][B][.65]{$2$}
		\psfrag{3}[B][B][.65]{$3$}
		\psfrag{a}[B][B][.65]{$\ge 2$}
		\psfrag{b}[B][B][.65]{+$1$}
		\psfrag{k}[B][B]{$k$}
		\psfrag{l}[B][B]{$l$}
		\psfrag{i}[B][B][.8]{\itshape inserting a vertex}
		\psfrag{s}[B][B][.8]{\itshape slicing}
	\centering\includegraphics[width=.95\linewidth]{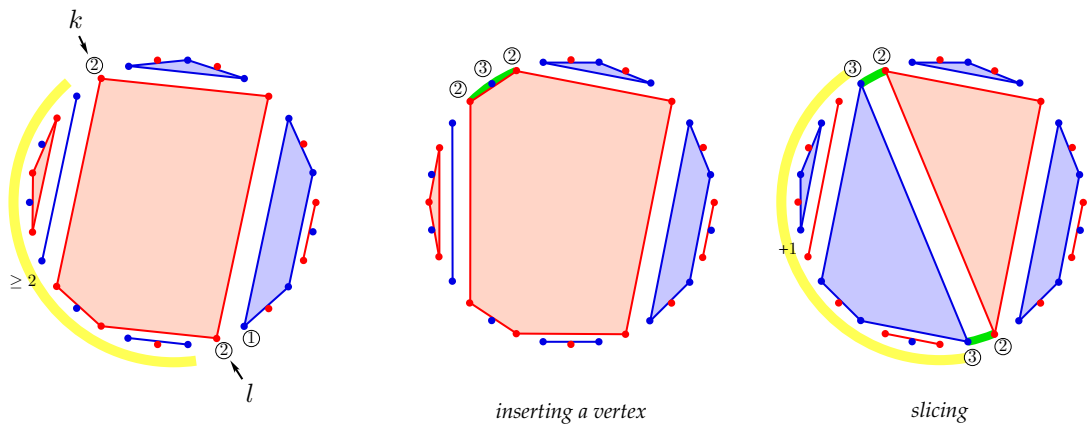}
	\caption{Consequences on the encoding Dyck paths of inserting a vertex or slicing at position~$k$.}
	\label{grow_l}
\end{figure}

These are exactly the moves considered in~\cite{Marchal03} (see in particular~\cite[Figure~1]{Marchal03}). The latter reference implies that, for each~$n$, the encoding Dyck path of~$\cP_n$ is uniformly distributed over the set of $2n$-step Dyck paths and that the convergence~\eqref{kaigh} holds almost surely for this choice of sequence $(\cP_n)_n$. As a result, there is no need to apply  Skorokhod's representation theorem in the proof of Theorem~\ref{thmpn}, so that~$\cP_{n}$ strongly converges toward the Brownian triangulation.
\end{pre}

\begin{pre}[Proof of Proposition~\ref{proppair}]
Inserting a short chord and inserting a long chord at position~$k$ respectively correspond on the unconstrained encoding Dyck path of Figure~\ref{encpair} to inserting one up-step followed by one down-step right after time~$k$ and lifting up by one the part of the path between time~$k$ and time~$l$. We conclude as above.
\end{pre}


\bibliographystyle{alpha}
\bibliography{main}

\end{document}